\documentclass{article}
\usepackage{amsmath,amsthm,amsopn,amstext,amscd,amsfonts,amssymb}
\usepackage{natbib}
\usepackage{verbatim}

\def\diag{\mathop{\rm diag}\nolimits}
\def\tr{\mathop{\rm tr}\nolimits}
\def\re{\mathop{\rm Re}\nolimits}

\def\etr{\mathop{\rm etr}\nolimits}

\newcommand{\Half}{\mbox{$\frac{1}{2}$}}

\renewenvironment{abstract}
                 {\vspace{6pt}
                  \begin{center}
                  \begin{minipage}{5in}
                  \centerline{\textbf{Abstract}}
                  \noindent\ignorespaces
                 }
                 {\end{minipage}\end{center}}
\newtheorem{thm}{\textbf{Theorem}}[section]
\newtheorem{cor}{\textbf{Corollary}}[section]
\newtheorem{lem}{\textbf{Lemma}}[section]
\theoremstyle{definition}
\newtheorem{defn}{\textbf{Definition}}[section]

\title{\huge \textbf{Compound and scale mixture of vector and spherical matrix variate elliptical distributions}}
\author{
  \textbf{Jos\'e A. D\'{\i}az-Garc\'{\i}a} \thanks{Corresponding author\newline
   {\bf Key words.}  Matricvariate, matrix variate, scale mixture, compound
    distribution, elliptical distribution, matrix variate hypergeometric distributions,
    matrix variate inverted hypergeometric distributions, zonal polynomials.\newline
    2000 Mathematical Subject Classification. 62E15, 60E05}\\
  Department of Statistics and Computation \\
  25350 Buenavista, Saltillo, Coahuila, Mexico \\
  E-mail: jadiaz@uaaan.mx \\[2ex]
  \textbf{Ram\'on Guti\'errez J\'aimez} \\
  Department of Statistics and O.R. \\
  University of Granada \\
  Granada 18071, Spain \\
  E-mail: rgjaimez@ugr.es\\
}
\date{}
\begin{document}
\maketitle
\begin{abstract}
Several matrix variate hypergeometric type distributions are
derived. The compound distributions of left-spherical matrix variate
elliptical distributions and inverted hypergeometric type
distributions with matrix arguments are then proposed. The scale
mixture of left-spherical matrix variate elliptical distributions
and univariate inverted hypergeometric type distributions is also
derived as a particular case of the compound distribution approach.
\end{abstract}

\section{\large{Introduction}}\label{Intro}

Four classes of matrix variate elliptical distributions have been defined and studied
by \citet{fz:90}. That $m \times n$ random matrix variate $\mathbf{X}$ is said to
have a matrix variate left-spherical distribution, the largest of the four classes
class of matrix variate elliptical distributions, if it density function is given by
$$
  \frac{c(m,n)}{|\mathbf{\Sigma}|^{n/2}|\mathbf{\Theta}|^{m/2}}h\left (\mathbf{\Sigma}^{-1/2}
  (\mathbf{X} - \boldsymbol{\mu})' \mathbf{\Theta}^{-1}(\mathbf{X} - \boldsymbol{\mu})
  \mathbf{\Sigma}^{-1/2}\right ),
$$
where $h$ is a real function, $c(m,n)$ denotes the normalization constant,
$\mathbf{\Sigma}$ is an $m \times m$ positive definite matrix, this fact being
denoted as $\mathbf{\Sigma} > \mathbf{0}$, $\mathbf{\Theta}$ is an $n \times n$
matrix, $\mathbf{\Theta} > \mathbf{0}$, and $\boldsymbol{\mu}$ is an $m \times n$
matrix.

When implementing frequentist or Bayesian methods, one may be interested in
situations where $\mathbf{X}$ has a density function of the following form, see
\citet{fz:90} and \citet{fl:99},
\begin{equation}\label{eliptica1}
    \frac{c(m,n)}{|\mathbf{\Sigma}|^{n/2}|\mathbf{\Theta}|^{m/2}}h\left (\mathbf{\Sigma}^{-1}
    (\mathbf{X} - \boldsymbol{\mu})' \mathbf{\Theta}^{-1}(\mathbf{X} - \boldsymbol{\mu})\right).
\end{equation}
This fact is denoted as $\mathbf{X} \sim \mathcal{ELS}_{m \times n}(\boldsymbol{\mu},
\mathbf{\Sigma}, \mathbf{\Theta}, h)$. This condition is equivalent to considering
the function $h$ as a symmetric function, i.e. $g: g(\mathbf{AB}) = g(\mathbf{BA})$
for any symmetric matrices $\mathbf{A}$ and $\mathbf{B}$. This condition is
equivalent to that in which $h(\mathbf{\mathbf{A}})$ depends on $\mathbf{A}$ only
through its eigenvalues, in which case the function $h(\mathbf{A})$ can be expressed
as $h(\lambda(\mathbf{A}))$, where $\lambda(\mathbf{A})= \diag(\lambda_{1}, \dots,
\lambda_{m})$ and $\lambda_{1}, \dots, \lambda_{m}$ are the eigenvalues of
$\mathbf{A}$. Two subclasses of matrix variate elliptical distributions are of
particular interest: the vector and spherical matrix variate elliptical
distributions. For these distributions, $\lambda (\mathbf{A}) \equiv \tr
(\mathbf{A})$ and $\lambda(\mathbf{A})$ represents  any function of eigenvalues of
$\mathbf{A}$, respectively. Note that vector matrix variate elliptical distributions
are a subclass of matrix variate spherical elliptical distributions. Many well-know
distributions are examples of these subclasses; one such is the matrix variate normal
distribution. Other variants include vector matrix variate elliptical distributions
e.g. Pearson type II, Pearson type VII, Kotz type , Bessel and Logistic, among many
others, see \citet{gv:93}. Yet other are matrix variate spherical elliptical
distributions, e.g. Pearson Type II, Pearson type VII and Kotz type among many
others, see \cite{fl:99}.

In the vectorial case, \citet[p. 33]{mh:82} proposed a means of generating a family
of vector variate elliptical distributions from a normal distribution. In the general
case, this idea has been extended to the matrix variate elliptical distributions by
\citet[pp. 78-79 and Section 4.1]{gv:93}. The situation in which the specific
elliptical distribution is a matrix variate normal is studied in \citet[Chapter
4]{gv:93}. Generically, distributions obtained by this procedure are termed scale
mixture matrix variate normal or elliptical distributions.

\citet{a:05} proposed the scale mixture of the vector variate Kotz type distribution,
also termed the $t$-type or generalised $t$ distribution. \citet{dggj:09} extend this
idea to matrix variate vector and spherical Kotz type distributions using two
approaches: scale mixture and compound matrix variate distributions.

These forms of obtaining vector or spherical matrix variate elliptical distributions
are of particular interest from a Bayesian point of view \citep{j:87,fl:99} and in
the context of shape theory, see \citet{cdg:08}.

This paper introduces several families of matrix variate elliptical distributions.
Section \ref{sec1} gives some results on integration, using zonal polynomials. In
terms of these results, various matrix variate hypergeometric type distributions are
proposed, as particular cases. These include well-known distributions such as central
and noncentral matrix variate inverted gamma (inverted Wishart) distributions and
matrix variate central and noncentral beta type II distributions. In Section
\ref{sec2}, assuming a hypergeometric type distribution for the matrix parameter in
matrix variate normal and matricvariate $T$ distributions, several families of matrix
variate elliptical distributions are obtained using the compound matrix variate
approach. Section \ref{sec3} introduces the scale mixture of a matrix variate
elliptical distribution, which is derived as a particular case of the compound matrix
variate approach, from where, all the results given in Section \ref{sec2} can be
particularised to this case.

\section{\large{Preliminary results}}\label{sec1}

Consider the following notation and definitions: the hypergeometric functions
${}_{p}F_{q}$ with matrix arguments are defined by
$$
  {}_{p}F_{q}(a_{1}, \dots, a_{p};b_{1}, \dots, b_{q}; \mathbf{Y}) = \sum_{k=0}^{\infty}\sum_{\kappa}
  \frac{(a_{1})_{\kappa} \cdots (a_{p})_{\kappa}}{(b_{1})_{\kappa} \cdots
  \mathbf{}  (b_{q})_{\kappa}} \ \frac{C_{\kappa}(\mathbf{Y})}{k!},
$$
where $\mathbf{Y}$ is a complex symmetric $m \times m$ matrix,
$C_{\kappa}(\mathbf{Y})$ is the zonal polynomial of $\mathbf{Y}$ of degree $k$,
$\kappa = (k_{1},\dots k_{m})$, $k_{1} \geq \cdots \geq k_{m} \geq 0$, $k_{1} +
\cdots + k_{m} = k$ and $a_{1}, \dots, a_{p},b_{1}, \dots, b_{q}$ are real or complex
constants,
$$
  (a)_{\kappa} = \prod_{i = 1}^{m}(a-(i-1)/2)_{k_{i}},
$$
with $(x)_{n} = x(x+1) \cdots (x+n-1)$, $(x)_{0} = 1$.

The multivariate gamma function is defined as
$$
  \Gamma_{m}[a] = \pi^{m(m-1)/4} \prod_{i = 1}^{m} \Gamma[a-(i-1)/2],
$$
and
$$
  \Gamma_{m}[a,\kappa] = \pi^{m(m-1)/4} \prod_{i = 1}^{m} \Gamma[a+k_{i}-(i-1)/2],
$$
where $\Gamma_{m}[a,\kappa] = (a)_{\kappa}\Gamma_{m}[a]$ with $\re(a) > (m-1)/2$, see
\citet{k:66} and \citet{mh:82}. Similarly,
\begin{eqnarray*}
  \Gamma_{m}[a,-\kappa] &=& \pi^{m(m-1)/4} \prod_{i = 1}^{m} \Gamma[a-k_{i}-(m-i)/2] \\
    &=& \pi^{m(m-1)/4} \prod_{i = 1}^{m} \Gamma[a-k_{m+1-i}-(i-1)/2],
\end{eqnarray*}
where $\re(a) > (m-1)/2 + k_{1}$. From
$$
    (-x)_{q} = (-1)^{q}(x-q+1)_{q} = \frac{(-1)^{q}\Gamma[x+1]}{\Gamma[x-q+1]},
$$
we obtain that
$$
  \Gamma_{m}[a,-\kappa] = \frac{(-1)^{k} \Gamma_{m}[a]}{(-a + (m+1)/2)_{\kappa}}.
$$

The multivariate beta function is defined as
$$
  \boldsymbol{\beta}_{m}[a,b] = \int_{\mathbf{0} < \mathbf{Y} < \mathbf{I}_{m}}|\mathbf{Y}|^{a-(m+1)/2}
  |\mathbf{I}_{m}-\mathbf{Y}|^{b-(m+1)/2}(d\mathbf{Y}) = \frac{\Gamma_{m}[a]\Gamma_{m}[b]}{\Gamma_{m}[a+b]},
$$
where $\re(a) > (m-1)/2$ and $\re(b) > (m-1)/2$, see \citet[p. 480]{h:55}.

Many distributions in multivariate analysis can be expressed in a form involving
hypergeometric functions with matrix arguments, as considered by \citet{r:75}, see
also \citet[Section 6.6.3, pp. 170-171]{p:82}. These distributions contain as
particular cases the central and noncentral gamma (Wishart) and matrix variate beta
type I and II distributions and are termed matrix variate hypergeometric gamma
(Wishart) type and matrix variate hypergeometric beta type I and II distributions. In
particular, for matrix variate hypergeometric beta type II distributions, we obtain
an alternative expression to the one given by \citet{r:75}, based on the following
lemma from \citet{k:66}.

\begin{lem}\label{lemma1}
If $\mathbf{R}$ is any arbitrary complex symmetric  $m \times m$ matrix, then
$$
  \int_{\mathbf{Y}>\mathbf{0}}|\mathbf{Y}|^{a -(m+1)/2}|\mathbf{I}_{m}+\mathbf{Y}|^{-(a+b)}
  C_{\kappa}(\mathbf{Y}\mathbf{R}) = \frac{(a)_{\kappa} \ \boldsymbol{\beta}_{m}[a,b]}{(-b +
  (m+1)/2)_{\kappa}}  C_{\kappa}(-\mathbf{R}),
$$
where $\re(a) > (m-1)/2 $ and $\re(b) > (m-1)/2 + k_{1}$.
\end{lem}
\textit{Proof.} As given in \citet{k:66}. \qed

If $\mathbf{Y} > \mathbf{0}$ has a matrix variate hypergeometric beta type II
distribution, the argument of the hypergeometric function involved, in Roux's
version, is $(\mathbf{I}_{m}+\mathbf{Y})^{-1}$, whereas in the version based on
Khatri's lemma, its argument is $\mathbf{Y}$, as we see below. The importance of this
fact is made apparent in the next section.

The next result is obtained immediately from Lemma \ref{lemma1}.

\begin{cor}\label{cor1}
Let $\mathbf{R}$ be any arbitrary complex symmetric  $m \times m$ matrix, then
$$
  \int_{\mathbf{Y}>\mathbf{0}}|\mathbf{Y}|^{a -(m+1)/2}|\mathbf{I}_{m}+\mathbf{Y}|^{-(a+b)}
  {}_{p}F_{q}(a_{1}, \dots, a_{p};b_{1}, \dots, b_{q}; \mathbf{R}\mathbf{Y})
  \hspace{2cm}
$$
$$\hspace{3cm}
  = \boldsymbol{\beta}_{m}[a,b]
  {}_{p+1}F_{q+1}(a_{1}, \dots, a_{p},a;b_{1}, \dots, b_{q}, -b +(m+1)/2; -\mathbf{R}),
$$
where $\re(a) > (m-1)/2$ and $\re(b) > (m-1)/2 + k_{1}$.
\end{cor}
\textit{Proof.} The final results are obtained using the fact that
$C_{\kappa}(b\mathbf{Y}) = b^{k}C_{\kappa}(\mathbf{Y})$, for a constant $b$. \qed

As a consequence of Corollary \ref{cor1} we have the following alternative definition
of the matrix variate hypergeometric beta type II distributions.

\begin{defn}\label{def1}
Let $\boldsymbol{\Xi}$ be any arbitrary complex symmetric  $m \times m$ matrix. Then
$\mathbf{Y}$ has a \textbf{matrix variate hypergeometric beta type II distribution}
if its density function is
$$
  f_{_{\mathbf{Y}}}(\mathbf{Y}) \propto |\mathbf{Y}|^{a -(m+1)/2}|\mathbf{I}_{m}+ \mathbf{Y}|^{-(a+b)}
  {}_{p}F_{q}(a_{1}, \dots, a_{p};b_{1}, \dots, b_{q}; \boldsymbol{\Xi} \mathbf{Y}),
  \quad \mathbf{Y} > \mathbf{0},
$$
where the constant of proportionality is
$$
   \frac{1}{\boldsymbol{\beta}_{m}[a,b] {}_{p+1}F_{q+1}(a_{1}, \dots, a_{p},a;b_{1}, \dots, b_{q}, -b +(m+1)/2;
   -\boldsymbol{\Xi})},
$$
with $\re(a) > (m-1)/2$ and $\re(b) > (m-1)/2 + k_{1}$.
\end{defn}

Let $\mathbf{Y}$ be a positive definite $m \times m$ matrix and let us define
$\mathbf{P} = \mathbf{Y}^{-1}$ then $(d\mathbf{Y}) =
|\mathbf{P}|^{-(m+1)}(d\mathbf{P})$. From \citet[eq. (4.1)]{r:75}, we have the
following result.

\begin{lem}\label{lemma2}\!
Let $\mathbf{\Upsilon}$ and $\boldsymbol{\Xi}$ be complex symmetric $m \times m$
matrices with $\re(\boldsymbol{\Xi}) > \mathbf{0}$. And assume that $\mathbf{Y}$ has
a hypergeometric matrix variate gamma type distribution. Then $\mathbf{P} =
\mathbf{Y}^{-1}$ has a \textbf{matrix variate inverted hypergeometric gamma type
distribution} with the following density function:
$$
  f_{_{\mathbf{P}}}(\mathbf{P}) \propto \etr\{-\boldsymbol{\Xi} \mathbf{P}^{-1}\}|\mathbf{P}|^{-a -(m+1)/2}
  {}_{p}F_{q}(a_{1}, \dots, a_{p};b_{1}, \dots, b_{q};\mathbf{\Upsilon} \mathbf{P}^{-1}),
  \quad \mathbf{P} > \mathbf{0},
$$
where the constant of proportionality is
$$
   \frac{|\boldsymbol{\Xi}|^{a}}{\Gamma_{m}[a] {}_{p+1}F_{q}(a_{1}, \dots, a_{p},a;b_{1}, \dots, b_{q};
   \mathbf{\Upsilon} \ \boldsymbol{\Xi}^{-1})},
$$
and $\re(a) > (m-1)/2$.
\end{lem}

Observe that if in Lemma \ref{lemma2}, $\mathbf{\Psi} = \mathbf{0}$, then
$\mathbf{P}$ has a  matrix variate central inverted gamma distribution. And if $p =
0$, $q = 1$ and $a = b_{1}$, then $\mathbf{P}$ has a  matrix variate noncentral
inverted Gamma distribution.

Similarly, from Definition \ref{def1}, we have

\begin{lem}\label{lemma3}
Let $\boldsymbol{\Xi}$ be any arbitrary complex symmetric  $m \times m$ matrix. Then
$\mathbf{P} = \mathbf{Y}^{-1}$ has a \textbf{matrix variate inverted hypergeometric
beta type II distribution} and its density function is
$$
  f_{_{\mathbf{P}}}(\mathbf{P}) \propto |\mathbf{P}|^{b -(m+1)/2}|\mathbf{I}_{m}+\mathbf{P}|^{-(a+b)}
  {}_{p}F_{q}(a_{1}, \dots, a_{p};b_{1}, \dots, b_{q};\boldsymbol{\Xi} \mathbf{P}^{-1}),
  \quad \mathbf{P} > \mathbf{0},
$$
where the constant of proportionality is
$$
   \frac{1}{\boldsymbol{\beta}_{m}[a,b] {}_{p+1}F_{q+1}(a_{1}, \dots, a_{p},a;b_{1}, \dots, b_{q}, -b +(m+1)/2;
   -\boldsymbol{\Xi})},
$$
and $\re(a) > (m-1)/2$ and $\re(b) > (m-1)/2 + k_{1}$.
\end{lem}
The distribution in Lemma \ref{lemma3} contains as particular cases the matrix
variate central and noncentral inverted beta type II distributions.

\citet{vr:74} studied other families of matrix variate hypergeometric distributions,
based on the ${}_{2}F_{1}(a,b;c;-\mathbf{Y})$ hypergeometric function with a matrix
argument. From this, and taking the limits of parameters $a$, $b$ or $c$, they obtain
the central and noncentral matrix variate gamma distributions, the central matrix
variate beta distribution and the matrix variate normal distribution. As we see in
\citet{ms:66} many other well-knows distributions can be obtained as particular cases
of this distribution. Next, we introduce this distribution using the multivariate
Mellin transform, \citet{m:97}.

Let $g(\mathbf{Y})$ be a function of the positive definite $m \times m$ matrix
$\mathbf{Y}$. The Mellin transform of $g(\mathbf{Y})$ is defined as
$$
  M_{g}(\mathbf{Y}) = \int_{\mathbf{Y} > \mathbf{0}} |\mathbf{Y}|^{\alpha-(m+1)/2}g(\mathbf{Y})
  (d\mathbf{Y}),
$$
where $\re(\alpha) > (m-1)/2$.
\begin{lem}\label{lemma4}
The Mellin transform of $g(\mathbf{Y}) = {}_{2}F_{1}(a,b;c;-\mathbf{Y})$ is given by
$$
  \int_{\mathbf{Y} > \mathbf{0}} |\mathbf{Y}|^{\alpha-(m+1)/2}{}_{2}F_{1}(a,b;c;-\mathbf{Y})
  (d\mathbf{Y}) = \frac{\boldsymbol{\beta}_{m}[\alpha,b-\alpha]\boldsymbol{\beta}_{m}
  [a-\alpha,c-\alpha]}{\boldsymbol{\beta}_{m}[a,c-a]},
$$
where $\re(\alpha) > (m-1)/2$, $\re(a-\alpha) > (m-1)/2$, $\re(b-\alpha) > (m-1)/2$,
$\re(c-\alpha)> (m-1)/2$ and $\re(c-a) > (m-1)/2$.
\end{lem}
\textit{Proof.} From the integral representation of ${}_{2}F_{1}$, see \citet[eq.
(2.12)]{h:55} and \citet[Theorem 7.4.2]{mh:82},
\begin{eqnarray*}
  M_{g}(\mathbf{Y}) &=& \int_{\mathbf{Y} > \mathbf{0}} |\mathbf{Y}|^{\alpha-(m+1)/2}
        {}_{2}F_{1}(a,b;c;-\mathbf{Y})(d\mathbf{Y}) \\
    &=& \displaystyle \frac{1}{\boldsymbol{\beta}_{m}[a,c-a]}\int_{\mathbf{Y} > \mathbf{0}}
        |\mathbf{Y}|^{\alpha-(m+1)/2}\int_{\mathbf{0} < \mathbf{R} < \mathbf{I}_{m}}
        |\mathbf{R}|^{a-(m+1)/2}\\
    & & \hspace{3cm}\times |\mathbf{I}_{m}-\mathbf{R}|^{c-a-(m+1)/2}|\mathbf{I}_{m}-
        \mathbf{Y}\mathbf{R}|^{-b}(d\mathbf{R})(d\mathbf{Y}) \\
    &=& \displaystyle \frac{1}{\boldsymbol{\beta}_{m}[a,c-a]}\int_{\mathbf{0} < \mathbf{R} < \mathbf{I}_{m}}
        |\mathbf{R}|^{a-(m+1)/2}|\mathbf{I}_{m}-\mathbf{R}|^{c-a-(m+1)/2}\\
    & & \hspace{3cm}\times \int_{\mathbf{Y} > \mathbf{0}} |\mathbf{Y}|^{\alpha-(m+1)/2}|\mathbf{I}_{m}-
        \mathbf{Y}\mathbf{R}|^{-b}(d\mathbf{Y})(d\mathbf{R}). \\
\end{eqnarray*}
Note that, $|\mathbf{I}_{m}-\mathbf{Y}\mathbf{R}| = |\mathbf{I}_{m}
-\mathbf{R}^{1/2}\mathbf{Y}\mathbf{R}^{1/2}| = |\mathbf{I}_{m}-\mathbf{W}|$, where
$\mathbf{R}^{1/2}$ is the positive definite square root of $\mathbf{R}$, such that
$\mathbf{R} = \left(\mathbf{R}^{1/2}\right)^{2}$ \citep[Theorem A9.3, p. 588]{mh:82}
and $|\mathbf{Y}| = |\mathbf{R}^{-1/2}\mathbf{W}\mathbf{R}^{-1/2}| =
|\mathbf{W}||\mathbf{R}|^{-1}$, with $\mathbf{W} =
\mathbf{R}^{1/2}\mathbf{Y}\mathbf{R}^{1/2}$. Then, $(d\mathbf{Y}) =
|\mathbf{R}|^{-(m+1)/2} (d\mathbf{W})$, from where
\begin{eqnarray*}
  M_{g}(\mathbf{Y}) &=& \displaystyle \frac{1}{\boldsymbol{\beta}_{m}[a,c-a]}
        \int_{\mathbf{0} < \mathbf{R} < \mathbf{I}_{m}}
        |\mathbf{R}|^{a-\alpha-(m+1)/2}|\mathbf{I}_{m}-\mathbf{R}|^{c-a-(m+1)/2}\\
    & & \hspace{3cm}\times \int_{\mathbf{W} > \mathbf{0}} |\mathbf{Y}|^{\alpha-(m+1)/2}
        |\mathbf{I}_{m}-\mathbf{W}|^{-b}(d\mathbf{W})(d\mathbf{Y}). \\
    &=& \displaystyle \frac{\boldsymbol{\beta}_{m}[\alpha,b-\alpha]}{\boldsymbol{\beta}_{m}[a,c-a]}
        \int_{\mathbf{0} < \mathbf{R} < \mathbf{I}_{m}} |\mathbf{R}|^{a-\alpha-(m+1)/2}
        |\mathbf{I}_{m}-\mathbf{R}|^{c-a-(m+1)/2}(d\mathbf{Y}) \\
    &=& \displaystyle \frac{\boldsymbol{\beta}_{m}[\alpha,b-\alpha]\boldsymbol{\beta}_{m}[a-\alpha,c-\alpha]}
    {\boldsymbol{\beta}_{m}[a,c-a]}.    \qed
\end{eqnarray*}

This was proved by \citet{vr:74}, who took the limit when $a$ tends to infinity in
Lemma \ref{lemma4} and found the Mellin transform of the function $g(\mathbf{Y}) =
{}_{1}F_{1}(b;c;-\mathbf{Y})$. Alternatively, the results can be obtained directly by
integration, as we shown below.

\begin{lem}\label{lemma5}
The Mellin transform of $g(\mathbf{Y}) = {}_{1}F_{1}(b;c;-\mathbf{Y})$ is given by
$$
  \int_{\mathbf{Y} > \mathbf{0}} |\mathbf{Y}|^{\alpha-(m+1)/2}{}_{1}F_{1}(b;c;-\mathbf{Y})(d\mathbf{Y}) =
  \frac{\Gamma_{m}[\alpha]\Gamma_{m}[c]\Gamma_{m}[b-\alpha]}{\Gamma_{m}[b]\Gamma_{m}[c-\alpha]},
$$
where $\re(\alpha) > (m-1)/2$, $\re(b-\alpha) > (m-1)/2$ and $\re(c-\alpha) >
(m-1)/2$.
\end{lem}
\textit{Proof.} Noting that from \citet[Theorem 7.4.2, p. 264]{mh:82}
\begin{eqnarray*}
\hspace{-0.75cm}
  {}_{2}F_{1}(c-b,\alpha;c;I_{m}) &=& \frac{1}{\boldsymbol{\beta}_{m}[a,c-a]}
        \int_{\mathbf{0} < \mathbf{R} < \mathbf{I}_{m}}
        |\mathbf{R}|^{c-b-(m+1)/2}|\mathbf{I}_{m}-\mathbf{R}|^{b-\alpha-(m+1)/2}(d\mathbf{R}) \\
    &=& \frac{\boldsymbol{\beta}_{m}[b-\alpha,c-b]}{\boldsymbol{\beta}_{m}[c-b,b]}.
\end{eqnarray*}
And from the Kummer relation discussed in \citet[Theorem 7.4.3, p. 265 and Theorem
7.3.4]{mh:82}, we have
\begin{eqnarray*}
  M_{g}(\mathbf{Y}) &=& \int_{\mathbf{Y} > \mathbf{0}} |\mathbf{Y}|^{\alpha-(m+1)/2}
        {}_{1}F_{1}(b;c;-\mathbf{Y})(d\mathbf{Y}) \\
    &=& \int_{\mathbf{Y} > \mathbf{0}} |\mathbf{Y}|^{\alpha-(m+1)/2}{}_{1} \etr(-\mathbf{Y})
        {}_{1}F_{1}(c-b;c;\mathbf{Y})(d\mathbf{Y}) \\
    &=& \Gamma_{m}[\alpha] {}_{2}F_{1}(c-b,\alpha;c;\mathbf{I}_{m})\\
    &=& \frac{\Gamma_{m}[\alpha]\Gamma_{m}[c]\Gamma_{m}[b-\alpha]}{\Gamma_{m}[b]\Gamma_{m}[c-\alpha]}. \qed
\end{eqnarray*}

Now, from Lemmas \ref{lemma4} and \ref{lemma5} taking $\mathbf{Y} =
\boldsymbol{\Xi}^{1/2} \mathbf{Y} \boldsymbol{\Xi}^{1/2}$ with $(d\mathbf{Y}) =
|\boldsymbol{\Xi}|^{(m+1)/2}(d\mathbf{Y})$, we obtain the following.

\begin{defn}\label{def2}
Let $\boldsymbol{\Xi}$ be any arbitrary complex symmetric  $m \times m$ matrix.
$\mathbf{Y}$ is said to have a \textbf{matrix variate generalised hypergeometric
distribution} if,
\begin{enumerate}
\item Its density function
is
$$
  f_{_{\mathbf{Y}}}(\mathbf{Y}) = \frac{|\boldsymbol{\Xi}|^{\alpha}\boldsymbol{\beta}_{m}[a,c-a]}
    {\boldsymbol{\beta}_{m}[\alpha,b-\alpha]\boldsymbol{\beta}_{m}[a-\alpha,c-\alpha]}
    |\mathbf{Y}|^{\alpha -(m+1)/2} {}_{2}F_{1}(a,b;c;-\boldsymbol{\Xi} \mathbf{Y}),
    \quad \mathbf{Y} > \mathbf{0},
$$
with $\re(\alpha) > (m-1)/2$, $\re(a-\alpha) > (m-1)/2$, $\re(b-\alpha) > (m-1)/2$,
$\re(c-\alpha) > (m-1)/2$ and $\re(c-a) > (m-1)/2$.
\item Or
$$
  f_{_{\mathbf{Y}}}(\mathbf{Y}) = \frac{|\boldsymbol{\Xi}|^{\alpha}\Gamma_{m}[b]\Gamma_{m}[c-\alpha]}
    {\Gamma_{m}[\alpha]\Gamma_{m}[c]\Gamma_{m}[b-\alpha]} |\mathbf{Y}|^{\alpha -(m+1)/2}
    {}_{1}F_{1}(b;c;-\boldsymbol{\Xi} \mathbf{Y}), \quad \mathbf{Y} > \mathbf{0},
$$
with $\re(\alpha) > (m-1)/2$, $\re(b-\alpha) > (m-1)/2$ and $\re(c-\alpha) >
(m-1)/2$.
\end{enumerate}
\end{defn}

\begin{lem}\label{lemma6}
Let $\boldsymbol{\Xi}$ be any arbitrary complex symmetric  $m \times m$ matrix. It is
said that $\mathbf{P} = \mathbf{Y}^{-1}$ has a \textbf{matrix variate inverted
generalised hypergeometric distribution} if,
\begin{enumerate}
\item Its density function
is
$$\hspace{-0.5cm}
  f_{_{\mathbf{P}}}(\mathbf{P}) = \frac{|\boldsymbol{\Xi}|^{\alpha}\boldsymbol{\beta}_{m}[a,c-a]}
    {\boldsymbol{\beta}_{m}[\alpha,b-\alpha]\boldsymbol{\beta}_{m}[a-\alpha,c-\alpha]}
    |\mathbf{P}|^{-\alpha -(m+1)/2} {}_{2}F_{1}(a,b;c;-\boldsymbol{\Xi} \mathbf{P}^{-1}),
    \quad \mathbf{P} > \mathbf{0},
$$
with $\re(\alpha) > (m-1)/2$, $\re(a-\alpha) > (m-1)/2$, $\re(b-\alpha) > (m-1)/2$,
$\re(c-\alpha) > (m-1)/2$ and $\re(c-a) > (m-1)/2$.
\item Or
$$
  f_{_{\mathbf{P}}}(\mathbf{P}) = \frac{|\boldsymbol{\Xi}|^{\alpha}\Gamma_{m}[b]\Gamma_{m}[c-\alpha]}
    {\Gamma_{m}[\alpha]\Gamma_{m}[c]\Gamma_{m}[b-\alpha]} |\mathbf{P}|^{-\alpha -(m+1)/2}
    {}_{1}F_{1}(b;c;-\boldsymbol{\Xi} \mathbf{P}^{-1}), \quad \mathbf{P} > \mathbf{0},
$$
with $\re(\alpha) > (m-1)/2$, $\re(b-\alpha) > (m-1)/2$ and $\re(c-\alpha) >
(m-1)/2$.
\end{enumerate}
\end{lem}
\textit{Proof.} Follows from Definition \ref{def2}, taking $\mathbf{P} =
\mathbf{Y}^{-1}$ with $(d\mathbf{Y}) = |\mathbf{P}|^{-(m+1)}(d\mathbf{P})$. \qed

\section{\large{Compound elliptical distribution of a random matrix}}\label{sec2}

In this section we propose several families of elliptical distributions based on an
extension to the matrix variate case of the vector case idea, introduced by
\citet{mh:82}, the approach Known as compound distribution. This approach was used by
\citet{r:71} and \citet{vr:74} for the distribution of a positive definite random
matrix.

In general, assume that the conditional distribution of
\begin{equation}\label{eq1}
  \mathbf{X}|\mathbf{P} \sim \mathcal{ELS}_{m \times n}(\boldsymbol{\mu},
  \mathbf{\Sigma}^{1/2} \mathbf{\Psi}(\mathbf{P})\mathbf{\Sigma}^{1/2},
  \mathbf{\Theta}, h),
\end{equation}
with $\mathbf{\Psi}:\Re^{m(m+1)/2}\rightarrow \Re^{m(m+1)/2}$,
$\mathbf{\Psi}(\mathbf{P}) > \mathbf{0}$; where $\mathbf{P}>\mathbf{0}$ has the
distribution function $G(\mathbf{P})$. Then $\mathbf{X}$ has a left-spherical
elliptical distribution (compound distribution) with a density function given by%
{\small
\begin{equation}\label{eliptica2}
    \frac{c(m,n)}{|\mathbf{\Sigma}|^{n/2}|\mathbf{\Theta}|^{m/2}}\int_{\mathbf{P}>\mathbf{0}}
    \frac{h\left (\mathbf{\Psi}(\mathbf{P})^{-1}\mathbf{\Sigma}^{-1/2}(\mathbf{X} - \boldsymbol{\mu})'
    \mathbf{\Theta}^{-1}(\mathbf{X} - \boldsymbol{\mu})\mathbf{\Sigma}^{-1/2}\right)
    dG(\mathbf{P})}{|\mathbf{\Psi}(\mathbf{P})|^{n/2}},
\end{equation}}
where $\mathbf{\Psi}(\mathbf{P})^{-1}$ denotes the inverse of the matrix
$\mathbf{\Psi}(\mathbf{P})$ (not the inverted function of $\mathbf{\Psi}(\cdot)$).

Let us now consider two particular matrix variate left-spherical elliptical
distributions, the matrix variate normal and the matricvariate $T$ distributions, see
\citet{di:67}, \citet[pp. 441-448]{bt:72} and \citet[pp. 138-141]{p:82}.

\section{\large{Compound matrix variate normal distribution}}

Recall that if $\mathbf{X} \sim \mathcal{N}_{m \times n}(\boldsymbol{\mu},
\mathbf{\Sigma}, \mathbf{\Theta})$, its density function is given by
$$
  \frac{1}{(2\pi)^{mn/2}|\mathbf{\Sigma}|^{n/2}|\mathbf{\Theta}|^{m/2}}\etr\left\{-\Half
  \mathbf{\Sigma}^{-1}(\mathbf{X} - \boldsymbol{\mu})' \mathbf{\Theta}^{-1}(\mathbf{X} -
  \boldsymbol{\mu})\right \}.
$$

\begin{thm}\label{teo1}
Assume that $\mathbf{X}|\mathbf{P}$ has a matrix variate normal distribution,
$$
  \mathbf{X}|\mathbf{P} \sim \mathcal{N}_{m \times n}(\boldsymbol{\mu},
  \mathbf{\Sigma}^{1/2}\mathbf{P}\mathbf{\Sigma}^{1/2},
  \mathbf{\Theta}),
$$
where $\mathbf{P}$ has a matrix variate inverted hypergeometric
gamma type distribution. By Lemma \ref{lemma2} its density function
is
$$
  g_{_{\mathbf{P}}}(\mathbf{P}) \propto \etr\{-\boldsymbol{\Xi} \mathbf{P}^{-1}\}|\mathbf{P}|^{-a -(m+1)/2}
  {}_{p}F_{q}(a_{1}, \dots, a_{p};b_{1}, \dots, b_{q};\mathbf{\Upsilon} \mathbf{P}^{-1}),
  \quad \mathbf{P} > \mathbf{0},
$$
where the constant of proportionality is
$$
   \frac{|\boldsymbol{\Xi}|^{a}}{\Gamma_{m}[a] {}_{p+1}F_{q}(a_{1}, \dots, a_{p},a;b_{1}, \dots, b_{q};
   \mathbf{\Upsilon} \ \boldsymbol{\Xi}^{-1})},
$$
and $\re(a) > (m-1)/2$. Then $\mathbf{X}$ has a matrix variate left-spherical
elliptical
distribution with density function%
{\small
$$
  \propto \frac{{}_{p+1}F_{q}\left(a_{1}, \dots, a_{p},a+\frac{n}{2};b_{1}, \dots, b_{q};
   \mathbf{\Upsilon}\left(\boldsymbol{\Xi} + \Half\mathbf{\Sigma}^{-1/2}(\mathbf{X} -
   \boldsymbol{\mu})'\mathbf{\Theta}^{-1}(\mathbf{X} - \boldsymbol{\mu})\mathbf{\Sigma}^{-1/2}
   \right)^{-1}\right)}{\left|\boldsymbol{\Xi} + \Half\mathbf{\Sigma}^{-1/2}(\mathbf{X} -
   \boldsymbol{\mu})'\mathbf{\Theta}^{-1}(\mathbf{X} - \boldsymbol{\mu})\mathbf{\Sigma}^{-1/2}\right |^{a+n/2}}
$$
with constant of proportionality
$$
   \frac{\Gamma_{m}[a + n/2]\ |\boldsymbol{\Xi}|^{a}}{(2\pi)^{mn/2}\Gamma_{m}[a] |\mathbf{\Sigma}|^{n/2}
   |\mathbf{\Theta}|^{m/2}
   {}_{p+1}F_{q}(a_{1}, \dots, a_{p},a;b_{1}, \dots, b_{q}; \mathbf{\Upsilon} \ \boldsymbol{\Xi}^{-1})}.
$$}
where $\re(a) > (m-1)/2$.
\end{thm}
\textit{Proof.} Follows immediately form \ref{eliptica2} and Lemma \ref{lemma2}. \qed

Observe that, by taking $\mathbf{\Upsilon} = \mathbf{0}$ in Theorem \ref{teo1} we
obtain that $\mathbf{X}$ has a matricvariate $T$ distribution, see \citet{di:67},
\citet[pp. 441-448]{bt:72} and \citet[pp. 138-141]{p:82}. Also, if we take $p=0$,
$q=1$  and $a = b_{1}$ we obtain that $\mathbf{X}$ has a noncentral matricvariate $T$
type 2 distribution. Then observing that
$$
  {}_{1}F_{1}(a;a; \mathbf{\Upsilon} \ \boldsymbol{\Xi}^{-1}) = \etr\{\mathbf{\Upsilon} \
  \boldsymbol{\Xi}^{-1}\}
$$
and by the Kummer relation \citep[eq. (6), p. 265]{mh:82}
$$
  {}_{1}F_{1}\left(a+ \frac{n}{2};a;
   \mathbf{\Upsilon}\left(\boldsymbol{\Xi} + \Half\mathbf{\Sigma}^{-1/2}(\mathbf{X} - \boldsymbol{\mu})'
   \mathbf{\Theta}^{-1}(\mathbf{X} - \boldsymbol{\mu})\mathbf{\Sigma}^{-1/2}\right)^{-1}\right)\hspace{3cm}
$$
$$
   = \etr\left\{\mathbf{\Upsilon}\left(\boldsymbol{\Xi} + \Half\mathbf{\Sigma}^{-1/2}
   (\mathbf{X} - \boldsymbol{\mu})'\mathbf{\Theta}^{-1}(\mathbf{X} -
   \boldsymbol{\mu})\mathbf{\Sigma}^{-1/2}\right)^{-1}\right\}\hspace{1cm}
$$
$$\hspace{2cm} \times
{}_{1}F_{1}\left(- \frac{n}{2};a;-
   \mathbf{\Upsilon}\left(\boldsymbol{\Xi} + \Half\mathbf{\Sigma}^{-1/2}(\mathbf{X} -
   \boldsymbol{\mu})'\mathbf{\Theta}^{-1}(\mathbf{X} - \boldsymbol{\mu})\mathbf{\Sigma}^{-1/2}
   \right)^{-1}\right).
$$
Observe that for $n/2$ an integer, ${}_{1}F_{1}$ is a polynomial of degree $mn/2$. In
this case the density of $\mathbf{X}$ is evaluated easily, see \citet[p. 258]{mh:82}.

From (\ref{eliptica2}) and Lemma \ref{lemma3} we have the following result.

\begin{thm}\label{teo2}
Assume that $\mathbf{X}|\mathbf{P} \sim \mathcal{N}_{m \times n}(\boldsymbol{\mu},
\mathbf{\Sigma}^{1/2}\mathbf{P}\mathbf{\Sigma}^{1/2}, \mathbf{\Theta})$, where
$\mathbf{P}$ has a matrix variate inverted hypergeometric beta type II distribution.
From Lemma \ref{lemma3}, its density function is
$$
  g_{_{\mathbf{P}}}(\mathbf{P}) \propto |\mathbf{P}|^{b -(m+1)/2}|\mathbf{I}_{m}+\mathbf{P}|^{-(a+b)}
  {}_{p}F_{q}(a_{1}, \dots, a_{p};b_{1}, \dots, b_{q};\boldsymbol{\Xi} \mathbf{P}^{-1}),
  \quad \mathbf{P} > \mathbf{0},
$$
where the constant of proportionality is
$$
   \frac{1}{\boldsymbol{\beta}_{m}[a,b] {}_{p+1}F_{q+1}(a_{1}, \dots, a_{p},a;b_{1}, \dots, b_{q}, -b +(m+1)/2;
   -\boldsymbol{\Xi})},
$$
and $\re(a) > (m-1)/2$ and $\re(b) > (m-1)/2 + k_{1}$. Then $\mathbf{X}$ has a matrix
variate left-spherical elliptical distribution with density function
$$
  \propto {}_{p+1}F_{q+1}\left(a_{1}, \dots, a_{p},a+n/2;b_{1}, \dots, b_{q}, -b +\frac{(m+1)}{2};
  \right. \hspace{3cm}
$$
$$\hspace{3cm}
  \left. \phantom{-b +\frac{(m+1)}{2}}-\boldsymbol{\Xi} + \Half\mathbf{\Sigma}^{-1/2}
  (\mathbf{X} - \boldsymbol{\mu})'\mathbf{\Theta}^{-1}(\mathbf{X} - \boldsymbol{\mu})
  \mathbf{\Sigma}^{-1/2}\right)
$$
where the constant of proportionality is
$$
   \frac{(2\pi)^{-mn/2}\boldsymbol{\beta}_{m}[a+n/2,b-n/2]|\mathbf{\Sigma}|^{-n/2}|\mathbf{\Theta}|^{-m/2}}
   {\boldsymbol{\beta}_{m}[a,b]
   {}_{p+1}F_{q+1}(a_{1}, \dots, a_{p},a+n/2;b_{1}, \dots, b_{q}, -b +(m+n+1)/2;-\boldsymbol{\Xi})},
$$
where $\re(a) > (m-1)/2$, $\re(b) > (m+n-1)/2 + k_{1}$.
\end{thm}

A result of particular interest is obtained from Theorem \ref{teo2} taking
$\boldsymbol{\Xi} = \mathbf{0}$. Similarly, in the Bayesian context, Theorem
\ref{teo2} generalises a result given in \citet{xu:90}, which can be obtained by
taking $\boldsymbol{\Xi} = \mathbf{0}$ and $p =q = 0$. In this latter case, by
applying the Kummer relation \citep[Theorem 7.4.3, p. 265]{mh:82} we obtain a
matricvariate confluent hypergeometric of the first kind distribution type.

\begin{thm}\label{teo3}
Assume that $\mathbf{X}|\mathbf{P} \sim \mathcal{N}_{m \times n}(\boldsymbol{\mu},
\mathbf{\Sigma}^{1/2}\mathbf{P}\mathbf{\Sigma}^{1/2}, \mathbf{\Theta})$, where
$\mathbf{P}$ has a matrix variate inverted generalised hypergeometric distribution.
By Lemma \ref{lemma6},
\begin{enumerate}
\item its density function is,%
$$\
  g_{_{\mathbf{P}}}(\mathbf{P}) \propto |\mathbf{P}|^{-\alpha -(m+1)/2} {}_{2}F_{1}(a,b;c;-\boldsymbol{\Xi}
  \mathbf{P}^{-1}), \quad \mathbf{P} > \mathbf{0},
$$
where the constant of proportionality is
$$
  \frac{|\boldsymbol{\Xi}|^{\alpha}\boldsymbol{\beta}_{m}[a,c-a]}{\boldsymbol{\beta}_{m}[\alpha,b-\alpha]
  \boldsymbol{\beta}_{m}[a-\alpha,c-\alpha]},
$$
with $\re(\alpha) > (m-1)/2$, $\re(a-\alpha) > (m-1)/2$, $\re(b-\alpha)
> (m-1)/2$, $\re(c-\alpha)
> (m-1)/2$ and $\re(c-a) > (m-1)/2$.
\item Or with density function
$$
  g_{_{\mathbf{P}}}(\mathbf{P}) \propto |\mathbf{P}|^{-\alpha -(m+1)/2} {}_{1}F_{1}(b;c;-\boldsymbol{\Xi}
  \mathbf{P}^{-1}), \quad \mathbf{P} > \mathbf{0},
$$
where the constant of proportionality is
$$
  \frac{|\boldsymbol{\Xi}|^{\alpha}\Gamma_{m}[b]\Gamma_{m}[c-\alpha]}{\Gamma_{m}[\alpha]
  \Gamma_{m}[c]\Gamma_{m}[b-\alpha]},
$$
with $\re(\alpha) > (m-1)/2$, $\re(b-\alpha) > (m-1)/2$ and $\re(c-\alpha) >
(m-1)/2$.
\end{enumerate}
Then $\mathbf{X}$ has a matrix variate left-spherical elliptical distribution and
\begin{enumerate}
\item its density function is
$$\
  \propto |\mathbf{\Sigma}^{-1}(\mathbf{X} - \boldsymbol{\mu})'\mathbf{\Theta}^{-1}
  (\mathbf{X} - \boldsymbol{\mu})|^{-(\alpha +n/2)}\hspace{6cm}
$$
$$\times
  {}_{3}F_{1}\left(a,b,\alpha +n/2 ;c;-2\boldsymbol{\Xi} \left(\mathbf{\Sigma}^{-1/2}
  (\mathbf{X} - \boldsymbol{\mu})'\mathbf{\Theta}^{-1}(\mathbf{X} - \boldsymbol{\mu})
  \mathbf{\Sigma}^{-1/2}\right)^{-1}\right),
$$
where the constant of proportionality is
$$
  \frac{2^{m\alpha}|\boldsymbol{\Xi}|^{\alpha}\Gamma_{m}[\alpha+n/2]\boldsymbol{\beta}_{m}[a,c-a]}
  {\pi^{mn/2}\boldsymbol{\beta}_{m}[\alpha,b-\alpha]\boldsymbol{\beta}_{m}[a-\alpha,c-\alpha]
  |\mathbf{\Sigma}|^{n/2}|\mathbf{\Theta}|^{m/2}},
$$
with $\re(\alpha) > (m-1)/2$, $\re(a-\alpha) > (m-1)/2$, $\re(b-\alpha)
> (m-1)/2$, $\re(c-\alpha)> (m-1)/2$ and $\re(c-a) > (m-1)/2$.
\item Or with density function given by
$$\
  \propto |\mathbf{\Sigma}^{-1}(\mathbf{X} - \boldsymbol{\mu})'\mathbf{\Theta}^{-1}
  (\mathbf{X} - \boldsymbol{\mu})|^{-(\alpha +n/2)}\hspace{6cm}
$$
$$\times
  {}_{2}F_{1}\left(a,\alpha +n/2 ;b;-2\boldsymbol{\Xi} \left(\mathbf{\Sigma}^{-1/2}
  (\mathbf{X} - \boldsymbol{\mu})'\mathbf{\Theta}^{-1}(\mathbf{X} - \boldsymbol{\mu})
  \mathbf{\Sigma}^{-1/2}\right)^{-1}\right),
$$
where the constant of proportionality is
$$
  \frac{2^{m\alpha}|\boldsymbol{\Xi}|^{\alpha}\Gamma_{m}[\alpha+n/2]|\boldsymbol{\Xi}|^{\alpha}
  \Gamma_{m}[b]\Gamma_{m}[c-\alpha]}{\pi^{mn/2}\Gamma_{m}[\alpha]\Gamma_{m}[c]\Gamma_{m}[b-\alpha]
  |\mathbf{\Sigma}|^{n/2}|\mathbf{\Theta}|^{m/2}},
$$
with $\re(\alpha) > (m-1)/2$, $\re(b-\alpha) > (m-1)/2$ and $\re(c-\alpha) >
(m-1)/2$.
\end{enumerate}
\end{thm}
\textit{Proof.} Follows from (\ref{eliptica2}) and Lemma \ref{lemma6}. \qed

\section{\large{Compound matricvariate $T$ distribution}}

From \citet{di:67}, \citet[pp. 441-448]{bt:72} and \citet[pp. 138-141]{p:82} we know
that $\mathbf{X}$ has a matricvariate $T$ distribution, denoting this fact as
$\mathbf{X} \sim \mathcal{MT}_{m \times n}(\nu, \boldsymbol{\mu}, \mathbf{\Sigma},
\mathbf{\Theta})$, if its density function is
$$
  \frac{\Gamma_{m}[(n+\nu)/2]}{\pi^{mn/2}\Gamma_{m}[\nu/2]|\mathbf{\Sigma}|^{n/2}|\mathbf{\Theta}|^{m/2}}
  \left|\mathbf{I}_{m} + \mathbf{\Sigma}^{-1}(\mathbf{X} - \boldsymbol{\mu})' \mathbf{\Theta}^{-1}
  (\mathbf{X} - \boldsymbol{\mu})\right|^{-(n+\nu)/2}.
$$
where $\nu > m-1$.
\begin{thm}\label{teo4}
Assume that $\mathbf{X}|\mathbf{P} \sim \mathcal{MT}_{m \times n}(\nu,
\boldsymbol{\mu}, \mathbf{\Sigma}^{1/2}\mathbf{P}\mathbf{\Sigma}^{1/2},
\mathbf{\Theta})$, where $\mathbf{P}$ has a matrix variate inverted hypergeometric
beta type II distribution with $\boldsymbol{\Xi} =\mathbf{0}$. From Lemma
\ref{lemma3}, its density function is
$$
  g_{_{\mathbf{P}}}(\mathbf{P}) \propto |\mathbf{P}|^{b -(m+1)/2}|\mathbf{I}_{m}+\mathbf{P}|^{-(a+b)},
  \quad \mathbf{P} > \mathbf{0},
$$
where the constant of proportionality is
$$
   \frac{1}{\boldsymbol{\beta}_{m}[a,b]},
$$
and $\re(a) > (m-1)/2$ and $\re(b) > (m-1)/2$. Then $\mathbf{X}$ has a matrix variate
left-spherical elliptical distribution with density function
$$
  \propto {}_{2}F_{1}\left(\frac{(n+\nu)}{2},a+\frac{n}{2}; -b +\frac{(m+n+1)}{2};\mathbf{\Sigma}^{-1}
  (\mathbf{X} - \boldsymbol{\mu})'\mathbf{\Theta}^{-1}(\mathbf{X} - \boldsymbol{\mu})\right)
$$
where the constant of proportionality is%
$$
   \frac{\Gamma_{m}[(n+\nu)/2]\boldsymbol{\beta}_{m}[a+n/2,b-n/2]}{\pi^{mn/2}\Gamma_{m}[\nu/2]
   \boldsymbol{\beta}_{m}[a,b] |\mathbf{\Sigma}|^{n/2}|\mathbf{\Theta}|^{m/2}
   },
$$
where $\re(a) > (m-1)/2$ and $\re(b) > (m+n-1)/2 + k_{1}$.
\end{thm}
\textit{Proof.} Follows from Lemma \ref{lemma3}, noting that
$$
  \left|\mathbf{I}_{m} + \mathbf{\Sigma}^{-1/2}(\mathbf{X} - \boldsymbol{\mu})'
  \mathbf{\Theta}^{-1}(\mathbf{X} - \boldsymbol{\mu})\mathbf{\Sigma}^{-1/2}
  \mathbf{P}^{-1}\right|^{-(n+\nu)/2} \hspace{4cm}
$$
$$ \hspace{2cm}
    = {}_{1}F_{0}\left(\frac{(n+\nu)}{2};-\mathbf{\Sigma}^{-1/2}(\mathbf{X} -
    \boldsymbol{\mu})'\mathbf{\Theta}^{-1}(\mathbf{X} - \boldsymbol{\mu})
    \mathbf{\Sigma}^{-1/2}\mathbf{P}^{-1}\right). \qed
$$
By applying the Euler relation \citep[eq. (7), p. 265]{mh:82} to results in Theorem
\ref{teo4}, the density function of $\mathbf{X}$ is then
$$
  \propto \left|\mathbf{\Sigma}^{-1}(\mathbf{X} - \boldsymbol{\mu})'\mathbf{\Theta}^{-1}
  (\mathbf{X} - \boldsymbol{\mu})\right|^{-(a+b+(n+\nu)/2-(m+1)/2)} {}_{2}F_{1}
  \left(-b-\frac{\nu}{2} + \frac{m+1}{2}\right.,
$$
$$
  \left.\phantom{-b-\frac{\nu}{2}\qquad}-a-b+\frac{m+1}{2}; -b +\frac{(m+n+1)}{2};
  \mathbf{\Sigma}^{-1}(\mathbf{X} - \boldsymbol{\mu})'\mathbf{\Theta}^{-1}(\mathbf{X}
  -\boldsymbol{\mu})\right).
$$
where, if $a$ and $b$ are integers, $2(a+b)> m+1$ and $m$ is odd, ${}_{2}F_{1}$ is a
polynomial of degree $m(a+b-(m+1)/2)$. Similarly, if $b$ and $\nu/2$ are integers,
$2(a+\nu/2)> m+1$ and $m$ is odd, ${}_{2}F_{1}$ is a polynomial of degree
$m(a+\nu/2-(m+1)/2)$, see \citet[p. 258]{mh:82}.

\section{\large{Scale mixture of elliptical distribution of a random matrix}}\label{sec3}

The approach known as the scale mixture of normal distributions, proposed by
\citet[p. 33]{mh:82} for the vector case and extended by \citet[Chapter 4]{gv:93} to
the matrix variate case, is obtained as a particular case of the approach described
in the Section \ref{sec2}. To do so, we take $m = 1$ in the distribution
$G(\mathbf{P})$, from where we obtain the following approach, termed the scale
mixture of an elliptical distribution, cited by \citet[pp. 78--79]{gv:93}.

Assume that the conditional distribution
\begin{equation}\label{eq2}
  \mathbf{X}|s \sim \mathcal{ELS}_{m \times n}(\boldsymbol{\mu},
  \phi(s)\mathbf{\Sigma}, \mathbf{\Theta}, h),
\end{equation}
where $\phi:(0, \infty)\rightarrow (0,\infty)$, with $s>0$ has the distribution
function $G(s)$. Then $\mathbf{X}$ has a left-spherical elliptical distribution
(scale mixture of elliptical distribution) with a density function given by
\begin{equation}\label{eliptica3}
    \frac{c(m,n)}{|\mathbf{\Sigma}|^{n/2}|\mathbf{\Theta}|^{m/2}}\int_{s>0}
    (\phi(s))^{-mn/2}h\left (\frac{1}{\phi(s)}\mathbf{\Sigma}^{-1}(\mathbf{X} - \boldsymbol{\mu})'
    \mathbf{\Theta}^{-1}(\mathbf{X} - \boldsymbol{\mu})\right)
    dG(s).
\end{equation}
As an example, consider the following version of Theorem \ref{lemma6} for $m = 1$.
\begin{thm}\label{teo5}
Assume that $\mathbf{X}|s$ has a matrix variate normal distribution,
$$
  \mathbf{X}|s \sim \mathcal{N}_{m \times n}(\boldsymbol{\mu}, s \mathbf{\Sigma},
  \mathbf{\Theta}),
$$
where $s$ has an inverted hypergeometric gamma type distribution. By Lemma
\ref{lemma2} its density function is
$$
  g_{_{S}}(s) \propto \exp\left\{-\frac{\xi}{s}\right\}s^{-a -1}
  {}_{p}F_{q}\left(a_{1}, \dots, a_{p};b_{1}, \dots, b_{q}; \frac{\upsilon}{s}\right),
  \quad s > 0,
$$
where $\upsilon > 0$, $\xi > 0$  and the constant of proportionality is
$$
   \frac{\xi^{a}}{\Gamma[a] {}_{p+1}F_{q}\left(a_{1}, \dots, a_{p},a;b_{1}, \dots, b_{q};
   \displaystyle\frac{\upsilon}{\xi}\right)},
$$
and $\re(a) > 0$. Then $\mathbf{X}$ has a matrix variate left-spherical elliptical
distribution with density function%
{\small
$$
  \propto \frac{{}_{p+1}F_{q}\left(a_{1}, \dots, a_{p},a+\displaystyle\frac{mn}{2};b_{1}, \dots, b_{q};
   \upsilon\left(\xi + \Half\tr\mathbf{\Sigma}^{-1}(\mathbf{X} -
   \boldsymbol{\mu})'\mathbf{\Theta}^{-1}(\mathbf{X} - \boldsymbol{\mu})
   \right)^{-1}\right)}{\left(\xi + \Half \tr\mathbf{\Sigma}^{-1}(\mathbf{X} -
   \boldsymbol{\mu})'\mathbf{\Theta}^{-1}(\mathbf{X} - \boldsymbol{\mu})\right )^{a+mn/2}}
$$
with constant of proportionality
$$
   \frac{\Gamma[a + mn/2]\ \xi^{a}}{(2\pi)^{mn/2}\Gamma[a] |\mathbf{\Sigma}|^{n/2}
   |\mathbf{\Theta}|^{m/2}
   {}_{p+1}F_{q}\left(a_{1}, \dots, a_{p},a;b_{1}, \dots, b_{q}; \displaystyle\frac{\upsilon}{\xi}\right)}.
$$}
where $\re(a) > 0$.
\end{thm}

Similarly, from Theorem \ref{teo5} particular cases  are obtained, taking, for
example, $\upsilon = 0$, but in this case we obtaining the matrix variate $T$
distribution (not the matricvariate $T$ distribution). Similar consequences are
obtained as particular cases, taking $m = 1$ from the Theorems \ref{teo2}-\ref{teo4}.

\section*{\large{Conclusions}}

This paper introduces several hypergeometric type distributions, which include many
that are well-known in the statistical literature, such as central and noncentral
matrix variate inverted Gamma (Wishart) and matrix variate beta type II
distributions, among many others.

Assuming that $\mathbf{P}$ has one of these hypergeometric distributions in
$$
  \mathbf{X}|\mathbf{P} \sim \mathcal{ELS}_{m \times n}(\boldsymbol{\mu},
  \mathbf{\Sigma}^{1/2} \mathbf{\Psi}(\mathbf{P})\mathbf{\Sigma}^{1/2},
  \mathbf{\Theta}, h),
$$
several left-spherical matrix variate elliptical families are
found. These enable us to study examples in which the tails of the distributions are
heavier or lighter than in the normal case. These approaches to obtaining elliptical
distributions are of particular interest from the Bayesian standpoint and for shape
theory distributions, see \citet{j:87} and \citet{cdg:08}, respectively.

\section*{\large{Acknowledgments}}

This research work was partially supported  by IDI-Spain, grants
FQM2006-2271 and MTM\-2008-05785, and CONACYT-M\'exico, research
grant no. \ 81512. This paper was written during J. A. D\'{\i}az-
Garc\'{\i}a's stay as a visiting professor at the Department of
Statistics and O. R. of the University of Granada, Spain.


\begin{thebibliography}{}

\bibitem[Arslan(2005)]{a:05}
    Arslan, O. 2005.
    A new class of multivariate distributions: Scale mixture of Kotz-type
    distributions.
    \textit{Statist. Prob. Letter}, \textbf{76}, 18-28.

\bibitem[Box and Tiao (1972)]{bt:72}
    Box, G. C. and Tiao, G. C. 1972.
    \textit{Bayesian Inference in Statistical Analysis.}
    Addison-Wesley, Reading, PA.

\bibitem[Caro-Lopera \textit{et al.} (2008)]{cdg:08}
    Caro-Lopera, F. J., D\'{\i}az-Garc\'{\i}a, J. A. and
    Gonz\'alez-Far\'{\i}as, G. 2008.
    Elliptical configuration distribution.
    \textit{J. Multivar. Ana.} To appear.

\bibitem[D\'{\i}az-Garc\'{\i}a and Guti\'errez-J\'aimez(2009)]{dggj:09}
    D\'{\i}az-Garc\'{\i}a, J. A.  and Guti\'errez-J\'aimez, R.
    2009.
    Compound and scale mixture of matricvariate and matrix variate Kotz-type
    distributions.
    \textit{J. Korean Statist. Soc.} To apear.

\bibitem[Dickey(1967)]{di:67}
    Dickey, J. M. 1967.
    Matricvariate generalizations of  the multivariate $t$- distribution and
    the inverted multivariate $t$-distribution.
    \textit{Ann. Math. Statist.} \textbf{38}, 511-518.

\bibitem[Fang and Li(1999)]{fl:99}
    Fang, K. T. and Li, R. 1999.
    Bayesian statistical inference on elliptical matrix
    distributions.
    \textit{J. Multivar. Anal.} \textbf{70}, 66-85.

\bibitem[Fang and Zhang(1990)]{fz:90}
    Fang, K. T. and Zhang, Y. T. 1990.
    \textit{Generalized Multivariate Analysis}.
    Science Press, Beijing, Springer-Verlang.

\bibitem[Herz(1955)]{h:55}
   Herz, C. S. 1955.
   Bessel functions of matrix argument.
   \textit{Ann. of Math}. \textbf{61}, 474-523.

\bibitem[Jammalamadaka \textit{et al.}(1987)]{j:87}
    Jammalamadaka, S. R., Tiwari, R. C. and Chib, S. 1987.
    Bayes prediction in the linear model with spherical symmetric
    errors.
    \textit{Econometrics Letters}, \textbf{24}, 39-44.

\bibitem[Gupta and Varga(1993)]{gv:93}
    Gupta, A. K. and Varga, T. 1993.
    \textit{Elliptically Contoured Models in Statistics}.
    Kluwer Academic Publishers, Dordrecht.

\bibitem[Khatri(1966)]{k:66}
   Khatri, C. G. 1966.
   On certain distribution problems based on positive definite quadratic
   functions in normal vector.
   \textit{Ann. Math. Statist}. \textbf{37}, 468-479.

\bibitem[Mathai(1997)]{m:97}
   Mathai, A. M. 1997.
   \textit{Jacobians of matrix transformations and functions of matrix
   argument}.
   World Scientific, London.

\bibitem[Mathai and Saxena(1966)]{ms:66}
   Mathai, A. M and Saxena, R. K. 1966.
   On a generalized hypegeometric distribution.
   \textit{Metrika} \textbf{11},  127--132.

\bibitem[Muirhead(1982)]{mh:82}
    Muirhead, R. J. 1982.
    \textit{Aspects of multivariate statistical theory}.
    Wiley Series in Probability and Mathematical Statistics,
    John Wiley \& Sons, Inc., New York.

\bibitem[Press (1982)]{p:82}
    Press, S. J. 1982.
    \textit{Applied Multivariate Analysis: Using Bayesian and Frequentist Methods of
    Inference}.
    Second Edition, Robert E. Krieger Publishing Company, Malabar, FL.

\bibitem[Roux(1971)]{r:71}
    J. J. J. Roux,
    On generalized multivariate distributions,
    S. Afr. Statist. J. 5 (1971) 91-100.

\bibitem[Roux(1975)]{r:75}
    Roux, J. J. J. 1975.
    \textit{New families of multivariate distributions.}
    In, G. P. Patil, S. Kotz, and J. K. Ord, (eds.) A Modern course on Statistical distributions
    in scientific work, Volume I, Model and structures, D. Reidel,
    Dordrecht-Holland, 281-297.

\bibitem[Van der Morwe and Roux(1974)]{vr:74}
    Van der Morwe, G. J. and Roux, J. J. J. 1974
    On generalized matrix-variate hypergeometric distribution.
    \textit{S. Afr. Statist. J}. \textbf{8}, 49-58.

\bibitem[Xu (1990)]{xu:90}
    Xu, J. L. 1990.
    Inverse Dirichlet distribution and its applications.
    In {\em Statistical Inference in Elliptically Contoured and Related
    Distributions}, (K. T. Fang, and T.~W. Anderson, Eds.) Allerton Press, New York,
    pp.~103--113.

\end{thebibliography}
\end{document}